 \documentclass[twoside,12pt]{amsart}
\usepackage{amsmath,latexsym,amssymb,times,mathptm,amsfonts,enumerate}
\usepackage{amsthm}

\def\proof{\noindent{\bf{Proof.} }}
\def\sqr#1#2{{\vcenter{\hrule height.#2pt
        \hbox{\vrule width.#2pt height#1pt \kern#1pt
                \vrule width.#2pt}
        \hrule height.#2pt}}}

\newtheorem{theorem}{Theorem}[section]

\newtheorem{lemma}[theorem]{Lemma}
\newtheorem{proposition}[theorem]{Proposition}

 \theoremstyle{definition}

\newtheorem{definition}[theorem]{Definition}

\newtheorem{remark}[theorem]{Remark}

\newtheorem{example}[theorem]{Example}

\newcommand{\f}[1]{\ensuremath{\mathfrak{#1}}}
\newcommand{\m}{\mathfrak{m}}

\newcommand{\height}{{\mbox{\rm{ht} }}}
\numberwithin{equation}{section}

\DeclareMathOperator{\depth}{\rm{depth \;}}
\newcommand{\core}[1]{\ensuremath{{\rm{core}}(#1)}}
\newcommand{\rees}[1]{\ensuremath{\mathcal{R}{(#1)}}}
\newcommand{\assg}[1]{\ensuremath{{\rm gr}_{#1}(R)}}
\newcommand{\spf}[1]{\ensuremath{\mathcal{F}{(#1)}}}

\begin{document}

\title{Reduction numbers and balanced ideals}
\author[L. Fouli]{Louiza Fouli}
\address{Department of Mathematical Sciences, New Mexico State University, Las Cruces, NM 88003, USA}
\email{lfouli@math.nmsu.edu}
\urladdr{http://www.math.nmsu.edu/~lfouli}

\subjclass[2010]{13A30, 13A15, 13B22}
\keywords{reductions, reduction number, core, associated graded ring, balanced}


\begin{abstract}
 Let $R$ be a Noetherian local ring and let $I$ be an ideal in $R$. The ideal $I$ is called balanced if the colon ideal $J:I$ is independent of the choice of the minimal reduction $J$ of $I$. Under suitable assumptions, Ulrich showed that $I$ is balanced if and only if the reduction number, $r(I)$, of $I$ is at most the `expected' one, namely $\ell(I)- \height I+1$, where $\ell(I)$ is the analytic spread of $I$. In this article we propose a generalization of balanced. We prove under suitable assumptions that if either $R$ is one-dimensional or the associated graded ring of $I$ is Cohen-Macaulay, then $J^{n+1}:I^n$ is independent of the choice of the minimal reduction $J$ of $I$ if and only if $r(I) \leq \ell(I)-\height I+n$.
 \end{abstract}

\maketitle

\section{Introduction}

Let $R$ be a Noetherian ring and let $I$ be an ideal in $R$. The {\it Rees algebra} $\rees{I}$ and  {\it the associated graded ring} $\assg{I}$ of $I$ are  $$\rees{I}=R[It]=\underset{i \geq 0} \oplus I^{i}t^{i} \mbox{ and } \assg{I}=R[It]/IR[It]=\underset{i \geq 0} \oplus I^{i}/I^{i+1}.$$
The projective spectrums of $\rees{I}$ and $\assg{I}$ are the blowup  of ${\rm Spec}(R)$ along $V(I)$ and the normal cone of $I$, respectively.  When studying various algebraic properties of these blowups a natural question to consider is which properties of the ring $R$ are transferred to these graded algebras. When $R$ is a local Cohen--Macaulay ring  and $I$ an ideal of positive height then if $\rees{I}$ is Cohen--Macaulay then so is $\assg{I}$, \cite{Hu1}. The converse does not hold true in general.  A celebrated theorem of Goto and Shimoda illustrates the intricate relationship between the Cohen--Macaulay property of these blowup algebras  and the reduction number of $I$. It states that when $(R, \m)$ is a local Cohen--Macaulay ring, with infinite residue field, dimension $d>0$ and $I$ an $\m$-primary ideal, then $\rees{I}$ is Cohen--Macaulay if and only if $\assg{I}$ is Cohen--Macaulay and the reduction number of $I$ is at most $d-1$, \cite{GoShi}. This theorem has inspired the work of many researchers and many generalizations  of it appeared in the literature in the late 1980s and early 1990s, see for example \cite{GHO},  \cite{HuHa1}, \cite{HuHa2}, \cite{GH}, \cite{JK}, \cite{AHT}, \cite{SUV}.

Recall that an ideal $J$ is a {\it reduction} of $I$ if $J \subset I$ and $\rees{I}$ is integral over $\rees{J}$ or equivalently if $J \subset I$
and $I^{n+1}=JI^{n}$ for some nonnegative integer $n$, see also Section~\ref{background}. The smallest non-negative integer such that the equality $I^{n+1}=JI^{n}$ holds is called the {\it reduction number of $I$ with respect to $J$} and is denoted by $r_J(I)$. 
When the ring is local then we consider  {\it minimal} reductions, where minimality is taken with respect to inclusion. In this case  the {\it reduction number of} $I$, denoted by $r(I)$,  is the minimum among all $r_J(I)$, where $J$ ranges over all minimal reductions of $I$.

An ideal $I$   satisfies the condition $ G_s$ for some integer $s$ if $\mu(I_{\f{p}}) \leq \dim R_{\f{p}}$ for every $\f{p} \in V(I)$ with $\dim R_{\f{p}} \leq
s-1$. The condition $G_s$  is  local and  rather mild. For example when $R$ is a Noetherian local ring with maximal ideal $\m$ and  dimension $d$, then any $\m$--primary ideal satisfies $G_{d}$. We say that an ideal $I$ satisfies $G_{\infty}$ if $I$ satisfies $G_s$ for every $s$.

Let $R$ be a local Gorenstein ring with  infinite residue field and  let $I$ be an ideal with $g= {\height} I >0$.  Suppose that $I$ satisfies $G_{\ell}$ and that  ${\depth}  R/I^j \geq \dim R/I -j+1$ for all $1 \leq j \leq \ell -g+1$, where $\ell=\ell(I)$ is the analytic spread of $I$.
In general, there are many classes of ideals that satisfy both the depth condition and $G_{\ell}$, for example ideals in the linkage class of a complete intersection satisfy these conditions; see \cite{CPU02} for more information.


A result of Johnson and Ulrich states that under the above conditions if
 $r(I) \leq \ell-g+ 1$ then $\assg{I}$ is Cohen--Macaulay. If in addition the height of $I$ is at least $2$ this also forces $\rees{I}$ to be Cohen--Macaulay,  \cite{JU}.  Moreover, the Castelnuovo--Mumford regularity of $\rees{I}$ and $\assg{I}$ can be calculated if $r(I) \leq \ell-g+1$.  The number $\ell(I)-\height I +1$ is known as the
 `expected reduction number' of $I$. This number was introduced by Ulrich in \cite{U2}, where he shows that under these assumptions an ideal $I$ has reduction number at most  the expected one if and only if the ideal is balanced. We say that $I$  is {\it balanced} if the colon ideal $J:I$ is independent of the minimal reduction $J$ of $I$, \cite[Theorem~4.8]{U2}. More precisely the definition of balanced is given below. 
  
 \begin{definition}{\rm (\cite[Definition~3.1]{U2}) Let $R$ be a Noetherian
local ring, let $I$ be an ideal, and let $s$ be a positive
integer. For a generating sequence $f_1, \ldots, f_n$ of $I$, let
$X$ be an $n$ by $n$ matrix of indeterminates, and write $[a_1,
\ldots, a_n]=[f_1, \ldots, f_n] \cdot X$ and $S=R(X)$. We say that
$I$ is {\it $s$--balanced} if there exist $n \geq s$ and $f_1,
\ldots, f_n$ as above such that $(a_{i_1}, \ldots, a_{i_s})S:IS$
yields the same $S$--ideal for every subset $\{i_1, \ldots, i_s \}
\subset \{1, \ldots, n\}$. }\end{definition}

We usually say that $I$ is balanced if $I$ is $\ell(I)$-balanced, where $\ell(I)$ is the analytic spread of $I$. It turns out that ideals that have the expected reduction number have many good properties. 
Next we discuss an application of \cite[Theorem~4.8]{U2}. Corso, Polini, and Ulrich make  use of the notion of balanced in order to establish a formula for the core of $I$. We recall here that \core{I}  is the intersection of all the reductions of $I$, see Section~\ref{background} for more details. Their theorem states that under the same assumptions as before one has that  $\core{I}=J(J:I)=J^2:I$ for all minimal reductions $J$ of $I$ if and only if $r(I) \leq \ell-g+1$. Therefore in this case the ideal $I$ is balanced if and only if $\core{I}=J^2:I$ for all minimal reductions $J$ of $I$. Most notably we see how the balanced condition, $J:I$ being independent of $J$, is intertwined with the formula for the core.

If the ideal is not balanced then it is natural to ask what can be a reasonable bound for the reduction number. The purpose of this article is to suggest a generalization of the notion of balanced and to establish bounds on the reduction number of an ideal in that case.



 We turn our attention to other known formulas for the core to obtain inspiration for what the generalization of balanced should be. Consider the following theorem due to Polini and Ulrich.

\begin{theorem}{\rm (\cite[Theorem~4.5]{PU})} \label{formulaPU} 
Let $R$ be a local Gorenstein ring with  infinite residue field $k$. Let $I$ be an ideal with $g= {\height} I >0$ and suppose that $I$ satisfies $G_{\ell}$ and that  ${\depth}  R/I^j \geq \dim R/I -j+1$ for all $1 \leq j \leq \ell -g$, where $\ell=\ell(I)$ is the analytic spread of $I$.
Let $J$ be a minimal reduction of $I$. If either ${\rm char}
 k =0$ or ${\rm char}  k > r_{J}(I)-\ell+g$ then  $\core{I}=J^{n+1}:I^{n} $ for  all   $n \geq \max \{r_{J}(I)-\ell+g,0\}$.
\end{theorem}

As one can see in Theorem~\ref{formulaPU} the characteristic of the residue field plays an important role when computing the core of an ideal. When appropriate we will be assuming  that the characteristic of the residue field is 0. In particular, under the set up of Theorem~\ref{formulaPU} the ideal $J^{n+1}:I^{n}$ is independent of the minimal reduction $J$ of $I$, since the formula for the core is independent of the choice of minimal reduction $J$ of $I$. Therefore, when
$n \geq {\max} \{r_{J}(I)-\ell+g,0\}$ then  $J^{n+1}:I^{n}$ is independent of the minimal reduction $J$ of $I$. Then it is natural to ask under which assumptions the converse holds true.

We propose the condition `$J^{n+1}:I^{n}$ is independent of the minimal reduction $J$' as a possible generalization of balanced. Notice that when $n=1$ then $I$ is balanced, by \cite[Theorem~2.6]{CPU02}.
The sequence of the ideals $\{J^{n+1}:I^{n}\}_{n \in \mathbb{N}}$ is decreasing as seen in Remark~\ref{decreasing seq}. 





We show that when the dimension of the ring $R$ is one then $J^{n+1}:I^{n}$ is independent of $J$ if and only if $n \geq r(I)$, Theorem~\ref{1dim balanced}. In the case of higher dimensions, we are able to show that the independence of the colon ideal $J^{n+1}:I^{n}$ from the choice of the minimal reduction $J$ of $I$ is equivalent to $r(I) \leq \ell(I)-\height I+n$, where $\ell(I)$ is the analytic spread of $I$, provided that $\assg{I}$ is Cohen-Macaulay, Theorem~\ref{BCM}.

\section{Background}\label{background}
Let $R$ be a Noetherian ring and $I$ an ideal in $R$. Recall that an deal $J$ is a {reduction} of $I$ if $J \subset I$ and $\rees{I}$ is integral over $\rees{J}$ or equivalently if $J \subset I$
and $I^{n+1}=JI^{n}$ for some nonnegative integer $n$. When the ring is local then we consider  { minimal} reductions, where minimality is taken with respect to inclusion. Northcott and Rees proved that if $R$ is a Noetherian local ring with maximal ideal $\m$ and infinite residue field then minimal reductions exist and either there are infinitely many or the ideal is basic, i.e. it is the only reduction of itself \cite{NR}. They show that minimal reductions correspond to Noether normalizations of the {\it special fiber ring of} $I$, $\spf{I}=\rees{I} \otimes R/\m$.

The concept of a reduction of an ideal was first introduced by Northcott and Rees in \cite{NR}, in order to facilitate the study of ideals and their powers.  
Reductions are in general smaller ideals with the same asymptotic behavior as the ideal $I$ itself. For example, all  minimal reductions  of $I$ have the same height and the same radical as $I$. Moreover, every minimal reduction $J$ of $I$ has the same minimal number of generators $\ell(I)$, where $\ell(I)$ is the {\it analytic spread} of $I$ and is defined to be the Krull dimension of the special fiber ring, $\spf{I}$,  of $I$.

Let $J$ be a minimal reduction of an ideal $I$ in a Noetherian local ring. The  \textit{reduction number of} $I$ \textit{with respect to} $J$, denoted by $r_J(I)$ is the smallest $n$ for which the equality $I^{n+1}=JI^{n}$ holds. This is denoted by
$r_{J}(I)$. In some sense the reduction number $r_{J}(I)$ measures how closely
related $J$ and $I$ are. The \textit{reduction
number} $r(I)$ of $I$ is the minimum of the reduction numbers
$r_{J}(I)$, where $J$ ranges over all minimal reductions of $I$.

In general, since an ideal has infinitely many reductions  it is natural to consider the {\it core} of the ideal, namely the intersection of all the (minimal) reductions of the ideal, \cite{RS}.
Several authors have determined formulas that describe the core in various settings, see for example  \cite{HS,  CPU01, CPU02, HySm1, HySm2,HT, PU, PUV, FPU1, FPU2}.

The core has many connections to geometry. For instance, Hyry and Smith have discovered a connection with a conjecture of  Kawamata on the non-vanishing of global sections of line bundles \cite{HySm1}. They prove that the validity of the conjecture is equivalent to a statement about core.

In a recent paper with  Polini and Ulrich we have uncovered yet another such connection with geometry. A scheme $X=\{P_{1}, \ldots, P_{s}\}$ of $s$ reduced points in $\mathbb{P}_{k}^{n}$ is said to have the {\it Cayley--Bacharach}
property if each subscheme of the form $X \backslash \{P_i\} \subset
\mathbb{P}_{k}^{n}$ has the same Hilbert function. It turns out that the structure of the core completely characterizes this property, namely $X$ has the Cayley--Bacharach property if and only if $\core{ \f{m}} =\m^{a+2}$, where $\m$ is the homogeneous maximal ideal of the homogeneous coordinate ring $R$ of $X$ and $a$ is the $a$--invariant of $R$, \cite{FPU2}.

We now discuss the notion of ideals of linear type. Recall that the Rees algebra $\rees{I}$ of $I$ is defined to be $\rees{I}=R[It]=\underset{i \geq 0} \oplus I^{i}t^{i}$. Let $R$ be a Noetherian ring and $I$ an ideal generated by $f_1, \ldots, f_n$. Then one has the following epimorphism $\phi: S=R[T_1, \ldots, T_n] \rightarrow \rees{I}$ given by $\phi(T_i)=f_it$. Let $J=\ker \phi$ and notice that $J$ is a graded ideal. Let $J=\oplus_{i=1}^{\infty} J_i$. Then $\rees{I}\simeq S/J$ and the ideal $J$ is often referred to as the defining ideal of $\rees{I}$. When $J=J_1$ then $I$ is called an ideal of {\it linear type}. It turns out that when $I$ is an ideal of linear type then $I$ is basic. The converse is not true in general.

The following is a well known result and we include it here for ease of reference.

\begin{lemma}\label{Jcancelation}
Let $R$ be a local Gorenstein ring and $I$ an ideal with 
$g={\height}  I>0$, $\ell=\ell(I)$,
and let $J$ be a minimal reduction of $I$. Assume that $I$ satisfies
$G_{\ell}$ and ${\depth}  R/I^{j} \geq \dim  R/I-j+1$
for $1\leq j \leq \ell -g$. Then
for every integer $n \geq 0$ and every integer $i \geq 0$
\[
J^{n+i}:J^{n}=J^{i}.
\]

\end{lemma}

\proof  First we note that ${\height}  J={\height}  I
>0$. According to \cite[Theorem~2.9]{U1} $I$ satisfies ${\rm
AN}_{\ell-1}$. Using $s=\ell-1$ in \cite[Theorem~1.11]{U1} we obtain
${\height}  J:I \geq \ell$ and hence $J$ satisfies $G_{\infty}$.
Hence by \cite[Remark~1.12]{U1} we have that $J$ satisfies ${\rm
AN}_{\ell-1}$. Using \cite[Theorem~1.8]{U1} we also obtain that $J$
satisfies sliding depth. Therefore ${\rm gr}_{J}(R)$ is
Cohen--Macaulay by \cite[Theorem~9.1]{HSV}. Then the cancelation is
clear, because ${\rm grade}  ({\rm gr}_{J}(R)_{+}) > 0$, since ${\rm
gr}_{J}(R)$ is Cohen--Macaulay and ${\height}  J >0$. \qed

We conclude this section with the following remark.

\begin{remark}\label{decreasing seq}{\rm
Let $R$ be a local Gorenstein ring and $I$ an ideal with 
$g={\height}  I>0$, $\ell=\ell(I)$,
and let $J$ be a minimal reduction of $I$. Assume that $I$ satisfies
$G_{\ell}$ and ${\depth}  R/I^{j} \geq \dim  R/I-j+1$
for $1\leq j \leq \ell -g$. Then $\{  J^{i+1}:I^{i} \}_{i\in
\mathbb{N}}$ is a decreasing sequence of ideals. To see this observe
that for all $i \geq 0$
\[
J^{i+1}:I^{i}=(J^{i+2}:J):I^{i}=J^{i+2}:JI^{i} \supset
J^{i+2}:I^{i+1},
\]
where the first equality holds according to Lemma~\ref{Jcancelation}.}
\end{remark}

\section{Main results}
We begin our investigation by considering the one-dimensional case. The first Lemma is analogous to
\cite[Lemma~4.7]{U2}.

\begin{lemma} \label{CF}
Let $R$ be an one--dimensional local Cohen--Macaulay ring with
canonical module $\omega_{R}$ and let $I$ be an ideal with
${\height}  I >0$. Assume that $I^i I^{-n} = a^i I^{-n}$ for some
$a \in I$ and for some positive integers $i$ and $n$, and that
$I^{r-1} \cong \omega_{R}$ for some positive integer $r$. Then
$I^{r+n} = a I^{r+n-1}$.
\end{lemma}

\proof First note that $a$ is a non--zerodivisor in $R$. Furthermore
we may assume $I^{r-1} = \omega_{R}$.  Since $I^i I^{-n} \supseteq
a^{i-1} I I^{-n} \supseteq a^i I^{-n}$ and $I^i I^{-n} = a^i I^{-n}$
we have that $a^{i-1} I I^{-n} = a^i I^{-n}$. Hence $I I^{-n} = a
I^{-n}$ since $a$ is a non--zerodivisor. Then for all $j>0$ it
follows that $I^j I^{-n} = a I^{j-1} I^{-n} = \cdots = a^j I^{-n}$.
For $j=r+n$ we obtain $a^{-r} I^{r+n} I^{-n} = a^n I^{-n}$ which
yields the following inclusions of fractional ideals:
\[\begin{array}{ll}
\hspace{1cm} a^{-r} I^{r+n} &\subset  a^{n} I^{-n}:I^{-n} \subset
R:I^{-n}=
(\omega_{R}: \omega_{R}): I^{-n}\\
& = \omega_{R} : (\omega_{R} I^{-n})
\stackrel{({\rm{*}})}{\subset} \omega_{R}: (a^{r-1} I^{-n})\\
&= a^{-r+1} \omega_{R}:(R:I^n) = a^{-r+1} \omega_{R}: ((\omega_{R} :
\omega_{R}):I^n)\\
& = a^{-r+1} \omega_{R}: (\omega_{R} : \omega_{R} I^n)
\stackrel{({\rm **})}{=} a^{-r+1} \omega_{R} I^n = a^{-r+1}
I^{r+n-1},
\end{array}\]
where $($*$)$ holds since $a^{r-1} \in \omega_{R}$ and $($**$)$
holds since $\dim R = 1$. Multiplication by $a^r$ implies that
$I^{r+n} \subset a I^{r+n-1}$ and thus $I^{r+n} = a I^{r+n-1}$. \qed

Using Lemma~\ref{CF} we are able to extend \cite[Theorem~2.6]{CPU02} in the case of a one--dimensional ring.

\begin{theorem}\label{1dim balanced}
Let $R$ be a one--dimensional local Gorenstein ring with residue
field of characteristic $0$. Let $I$ be an ideal with
${\height} I >0$ and $J$ a minimal reduction of $I$.
Then the following are equivalent for a positive integer $n$ $:$
\begin{enumerate}[$($a$)$]
\item $J^{n+1}:I^{n}$ is independent of $J$;
\item $ \core{I}=J^{n+1}:I^{n}$ for some $J$;
\item $n \geq r(I)$.
\end{enumerate}
\end{theorem}

\proof Notice that ${\height} {I}=1=\ell(I)$ since $\dim R=1$. By
\cite[Theorem~2.1]{Huc} we have that $r_{J}(I)$ is independent of the
minimal reduction $J$ of $I$. Hence $r_{J}(I)=r(I)$. Let
$r=r_J(I)=r(I)$. 

Suppose that $n \geq r$. Then by \cite[Theorem~4.5, Remark~4.8]{PU} we have that 
$\core{I}=J^{n+1}:I^{n}$ and $J^{n+1}:I^{n}$ is
independent of the minimal reduction $J$ of $I$. This establishes the implications $($c$) \Rightarrow ($a$)$ and $($c$) \Rightarrow ($b$)$.

To prove (b) $\Rightarrow$ (c) suppose that $\core{I}=J^{n+1}:I^{n}$.
By \cite[Theorem~4.5]{PU} we know that $\core{I}=J^{m+1}:I^{m}$ for
$m \geq r$. Let $m \geq \max \{r,n\}$. Then
\begin{eqnarray*}
\core{I}&=&J^{m+1}:I^{m} \subset
J^{m+1}:J^{m-n}I^{n}\\
&=&(J^{m+1}:J^{m-n}):I^{n}
\stackrel{(1)}{=}J^{n+1}:I^{n}=\core{I},
\end{eqnarray*}
where (1) holds since $J$ is generated by a single regular element.
Therefore $J^{m+1}:I^{m}=J^{m+1}:J^{m-n}I^{n}$. Since $R$ is
Gorenstein then by linkage we have $I^{m}=J^{m-n}I^{n}$. Hence $n \geq r$.

Finally, in order to prove that (a) $\Rightarrow$ (c)  notice that there exists $m \gg 0$ such that for
general linear combinations $f_{1}, \ldots, f_{m}$ of the generators
of $I$, we have that $(f_{i})$ forms a reduction of $I$ for $1 \leq
i \leq m$ and $I^{n+1}=(f_{1}^{n+1}, \ldots, f_{m}^{n+1})$ since
${\rm char}  k=0$. For example one may take $m=e(R)$, the
multiplicity of the ring $R$. Let $J=(a)$. Then for all $1 \leq i
\leq m$,
\[
a^{n+1}I^{-n}=a^{n+1}:I^{n}=f_{i}^{n+1}:I^{n}= f_{i}^{n+1}I^{-n}.
\]Hence $a^{n+1}I^{-n}=I^{n+1}I^{-n}$. Then by Lemma \ref{CF} we
obtain $I^{n+1}=aI^{n}$ and  thus $n \geq r$. \qed

Next we give a description for the canonical module of the extended Rees ring.

\begin{remark}\label{omegaB}
{\rm Let $R$ be a local Gorenstein ring and $I$ an ideal with
$g={\height}  I>0$, $\ell=\ell(I)$, and $J$ a minimal reduction of
$I$. Write $B=R[It,t^{-1}]$. Assume that $I$ satisfies $G_{\ell}$
and ${\depth}  R/I^{j} \geq \dim  R/I-j+1$ for $1\leq
j \leq \ell -g$. We fix a graded canonical module for the ring $B$
such that $\omega_{B} \subset R[t,t^{-1}]$ and
$[\omega_{B}]_{i}=Rt^{i}$ for all $i \ll 0$. Notice that this
uniquely determines $\omega_{B}$ as a submodule of $R[t,t^{-1}]$.
According to \cite[Remark~2.2]{PU} we have a description of $\omega_{B}$. For all $n \geq \max
\{r_J(I)-\ell+g, 0\}$}
\[
\omega_{B}= \bigoplus_{i\in \mathbb{Z}}  (J^{i}:_{R}
I^{n})t^{i-n+g-1}= \ldots \oplus Rt^{g-n-1} \oplus (J:I^{n})t^{g-n}
\oplus \ldots .
\]
\end{remark}

Let $R$ be a Noetherian local ring that is an epimorphic image of a local Gorenstein ring.
Let $B$ be a $\mathbb{Z}$--graded Noetherian $R$--algebra with
$B_{0}=R$ and unique homogeneous maximal ideal $\m$. We also
assume that $B/\m$ is a field. Let $\omega_{B}$ be the graded
canonical module of $B$. Recall that  the {\it a--invariant} of $B$ is
$a(B)= - {\rm indeg}(\omega_{B} \otimes_{B} B/\m)$. Notice that if $B$ is positively graded then $a(B)=-{\rm indeg} \omega_{B}$. 

In the setting of Theorem~\ref{1dim balanced} the reduction numbers were independent of the choice of minimal reduction as seen in the proof of Theorem~\ref{1dim balanced}.  In the next Proposition we provide conditions that guarantee the independence of the reduction numbers. This result was known in the case $I$ is equimultiple and $\depth \assg{I}_{+}\geq \dim R-1$, \cite[Theorem~2.1]{Huc}. In the case that $I$ is an $\m$--primary ideal this result was also obtained by \cite[Theorem~1.2]{Tr}. Our setup is  more general.

\begin{proposition}\label{r=rJ}
Let $R$ be a local Gorenstein ring with infinite residue field. Let
$I$ be an ideal with $g= \height{I} >0$ and $\ell=\ell(I)$.
Assume that $I$ satisfies
$G_{\ell}$ and ${\depth}  R/I^{j} \geq \dim  R/I-j+1$
for $1\leq j \leq \ell -g$. We further assume that ${\rm gr}_{I}(R)$
is Cohen--Macaulay. Then
$r(I)=r_{J}(I)$ for every minimal reduction $J$ of $I$. 

\end{proposition}

\proof According to \cite[Corollary~5.5]{JU} either $r(I)=0$ or $r(I) >
\ell-g$. If $r(I)=0$ then there is nothing to show. So we assume
that $r(I)
>\ell-g$. Let $J$ be a minimal reduction of $I$. Then $r_{J}(I) > \ell-g$.

Let $\f{p} \in {\rm Spec} (R)$ such that $\f{p} \supset I$ and
$\ell(I_{p})={\height}  \f{p} < \ell$. Then $I_{\f{p}}$ is of linear type and thus
$r(I_{\f{p}})=0$, according to \cite[Proposition~1.11]{U1}. Thus $r(I_{\f{p}})- {\height}  \f{p} \leq -g <
r_{J}(I)-\ell$ and hence $a({\rm gr}_{I}(R))=r_{J}(I)-\ell$ by
\cite[Corollary~4.5]{AHT}. But the $a$-invariant of $\assg{I}$ is independent of the choice of the minimal reduction $J$ and thus $r_{J}(I)$ is independent of
$J$.  Hence $r_{J}(I)=r(I)$. \qed

The next result is essentially obtained in \cite[Corollary~2.4]{U2} but we are able to weaken the assumptions on the depth condition.

\begin{proposition}\label{mBCM} 
Let $R$ be a local Gorenstein ring with residue field of
characteristic $0$. Let $I$ be an ideal with $g= {\height} {I}
>0$, $\ell=\ell(I)$, and let $J$ be a minimal reduction of $I$. Suppose that
$I$ satisfies $G_{\ell}$ and $\depth  R/I^{j} \geq
\dim  R/I-j+1$ for $1\leq j \leq \ell -g$. We further assume
that ${\rm gr}_{I}(R)$ is Cohen--Macaulay. Then
\begin{enumerate} [$($a$)$]
\item $J:I^{n} \neq R \mbox{ for all } n \leq \max\{r(I)-\ell+g, 0\}$  and

\item $\max \{r(I)-\ell+g, 0\}=\min \{ i \mid I^{i+1} \subset
\core{I} \}$. 
\end{enumerate}
\end{proposition}

\proof Write $G={\rm
gr}_{I}(R)$ and $B=R[It,t^{-1}]$. As $G$ is Cohen--Macaulay then so
is $B$, since $G={\rm gr}_{I}(R) \simeq B/(t^{-1})$. According to
Proposition~\ref{r=rJ} one has that $r_{J}(I)=r(I)$. Furthermore
$a(G)=\max \{r(I)-\ell,-g \}$ by \cite[Theorem~3.5]{SUV}. On the
other hand, $a(G)= a(B)-1$ since $G$ is Cohen--Macaulay and $G \simeq
B/(t^{-1})$. Therefore $a(B)=m-g+1$, where $m=\max \{ r(I)-\ell(I)+g, 0\}$. Hence
$[\omega_{B}]_{g-m-1}=R$ and $J:I^{m}=[\omega_{B}]_{g-m} \neq R$ by
Remark~\ref{omegaB}, since $r_{J}(I)=r(I)$.  Hence $J:I^{n} \subset
J:I^{m} \neq R$ for all $n \leq m$. This proves $($a$)$.

For part $($b$)$ we claim that
$m=\min \{i\mid I^{i+1} \subset \core{I}\}$.
To see this observe that $J \subset J:I^{m}$ and $J:I^{m}$ is
independent of $J$ by \cite[Remark~2.3]{PU}, since
$r_{J}(I)=r(I)$. Thus $I \subset J:I^{m}$ and hence $I^{m+1} \subset
J$. Consequently $I^{m+1} \subset \core{I}$. But since $J:I^{m} \neq
R$ we have that $I^{m} \not\subset J$ and therefore $I^{m}
\not\subset \core{I}$. \qed

In order to extend Theorem~\ref{1dim balanced} we prove the first two statements are equivalent in higher dimensions without any additional assumptions on the associated graded ring of the ideal. 

\begin{proposition}\label{coreiffindep}
Let $R$ be a local Gorenstein ring with residue field of
characteristic $0$. Let $I$ be an ideal with $g= {\height} {I}
>0$, $\ell=\ell(I)$, and let $J$ be a minimal reduction of $I$. Suppose
that $I$ satisfies $G_{\ell}$ and $\depth  R/I^{j} \geq
\dim  R/I-j+1$ for $1\leq j \leq \ell -g$. Then the
following are equivalent for an integer $n$ $:$
\begin{enumerate}[$($a$)$]
\item $J^{n+1}:I^{n}$ is independent of  $J$;
\item $\core{I}=J^{n+1}:I^{n}$ for every $J$.
\end{enumerate}
\end{proposition}

\proof By \cite[Theorem~4.5]{PU} we have that
$\core{I}=J^{m+1}:I^{m}$ for $m \gg 0$ and any minimal reduction $J$
of $I$. 

Suppose that $J^{n+1}:I^{n}$ is independent of  $J$.
Notice that $J^{n+1}:I^{n}
\subset J^{n+1}:J^{n}=J$, where the equality holds by Lemma~\ref{Jcancelation}.
Since $J^{n+1}:I^{n}$ is independent
of $J$ it follows that $J^{n+1}:I^{n} \subset
\core{I}=J^{m+1}:I^{m}$ for $m \gg 0$. By Remark~\ref{decreasing seq} we have that $\{  J^{i+1}:I^{i}
\}_{i\in \mathbb{N}}$ is a decreasing sequence of ideals and hence it
follows that $\core{I}=J^{n+1}:I^{n}$ for every minimal reduction $J$ of $I$.

The other implication is clear since the formula for $\core{I}$
is independent of the choice of the  minimal reduction $J$ of $I$. \qed

We are now ready to prove the main result of this article. If we assume that the associated graded ring of $I$ is Cohen--Macaulay
then we obtain a generalization to Theorem~\ref{1dim balanced} in higher dimensions.

\begin{theorem}\label{BCM}
Let $R$ be a local Gorenstein ring with residue field of
characteristic $0$. Let $I$ be an ideal with $g= {\height} {I}
>0$, $\ell=\ell(I)$, and let $J$ be a minimal reduction of $I$. Suppose
that $I$ satisfies $G_{\ell}$ and $\depth  R/I^{j} \geq
\dim  R/I-j+1$ for $1\leq j \leq \ell -g$. We further assume
that ${\rm gr}_{I}(R)$ is Cohen--Macaulay. Then the following are
equivalent for an integer $n$ $:$
\begin{enumerate}[$($a$)$]
\item $J^{n+1}:I^{n}$ is independent of  $J$;
\item $ \core{I}=J^{n+1}:I^{n}$ for every $J$;
\item $n \geq {\max} \{r(I)-\ell+g, 0\}$.
\end{enumerate}
\end{theorem}

\proof
The first two statements are equivalent by Proposition~\ref{coreiffindep}.
 Write $G={\rm gr}_{I}(R)$
and $B=R[It,t^{-1}]$. Since $G$ is Cohen--Macaulay then so is $B$
since $G={\rm gr}_{I}(R) \simeq B/(t^{-1})$. Notice that
$r_{J}(I)=r(I)$ by Proposition~\ref{r=rJ}. 

Let $m=\max \{ r(I)-\ell+g, 0\}$ and suppose that $n \geq m$.  Then
$\core{I}=J^{n+1}:I^{n}$, for any minimal reduction $J$ of $I$ according to \cite[Theorem~4.5]{PU}, since $r_{J}(I)=r(I)$. 

Finally, suppose that $\core{I}=J^{n+1}:I^{n}$. Then
$J^{n+1} \subset \core{I}$ for every minimal reduction $J$ of $I$.
Since ${\rm char} k =0$ we obtain that $I^{n+1} \subset \core{I}$.
Therefore $n \geq m$ by Proposition~\ref{mBCM}.\qed

The following example is due to Angela Kohlhass. It establishes that without the Cohen-Macaulay assumption on the associated graded ring the result of Theorem~\ref{BCM} does not hold in general.

\begin{example}[A. Kohlhass]\label{Angela Ex}{\rm
Let $R=k[[x,y]]$ be a power series ring over a field $k$ of characteristic $0$. Let $I=(x^{10},x^4y^5,y^9)$ and $J$ a general  minimal reduction of $I$. Then  $I$ is $\m$-primary, where $\m$ is the maximal ideal of $R$, $r(I)=4$, and  $\depth \assg{I}=0$. 
It turns out that  $J^4:I^3  = \core{I}= J^5:I^4$. 

}
\end{example}

\begin{remark}{\rm
We remark that in Example~\ref{Angela Ex} the associated graded ring of the ideal $I$ has depth 0 and the ideal $J^4:I^3$ is independent of the choice of the minimal reduction $J$ of $I$, whereas $r(I)=4$. This shows that in general Theorem~\ref{BCM} does not hold without any assumptions on $\assg{I}$. It is conceivable that when $\depth \assg{I} \geq \dim R-1$ then a similar statement might hold.
}
\end{remark}


\end{document}